\documentclass{article}
\usepackage{amssymb}
\usepackage{amstext}
\usepackage{amsthm}
\usepackage{a4wide}
\usepackage{amsmath}
\newcommand{\su}{\sum_{i\in I} a_i}
\newcommand{\F}{{\mathfrak F}}
\renewcommand{\O}{{\mathfrak O}}
\newcommand{\nn}{{\mathbb N}^{(\infty)}}
\newcommand{\nnq}{{\mathbb N}^{(\infty)}\langle\!\langle S^* \rangle\!\rangle}
\newcommand{\nnp}{{\mathbb N}\langle S ^* \rangle}

\def\xx#1 {\newtheorem{#1}[thm]{#1}}
\xx Theorem
\xx Fact
\xx Definition
\xx Corollary
\xx Remark 
\xx Example

\begin{document}

\author{Martin Goldstern}

\title{Completion of semirings\footnote{
In 1985 I wrote a diploma thesis (in German) on ``completion of semirings'' 
at the Institute of Algebra at the Vienna University of Technology. 
My advisor was Werner Kuich.
\endgraf
I never published my thesis, but 
 some people have expressed an 
interest in seeing it, and references to it appear in a few published 
papers or books  (\cite{Golan:1992}, \cite{Kuich:1991},\cite{Karner:1993},
\cite{Kuich:1997}, \cite{Glazek:2002}).     
  So here is an English summary.}}

\date{Spring 1985; Summer 2002}

\maketitle

\nocite{Golan+Wang:1996}
\nocite{Golan:1999}

\nocite{Kuich:1997}
\nocite {Kuich:1995}
\nocite{Kuich:1991}

\begin{abstract}
A semiring can be ``completed''
(i.e., embedded into a semiring in which all infinite sums are 
defined and satisfy some reasonable properties) iff this semiring can 
be naturally partially ordered.       This construction is ``natural'' 
(a left adjoint to the forgetful functor), and quite straightforward.
\end{abstract}

\section*{Definitions and easy facts}

\begin{Definition}
A semiring is a structure $(S, {+}, {\cdot}, 0,1)$ such that
 $(S,{+},0)$ is a commutative monoid, $(S,{\cdot},1)$ is a monoid, 
and the two distributive laws $x\cdot (y+z)=x\cdot y + x\cdot z$, 
$(y+z)\cdot  x = y \cdot  x +  z\cdot  x $ hold.

\end{Definition}
\begin{Definition}
A complete semiring
 $(S, {+}, {\cdot}, 0,1, \Sigma  )$ is a semiring in which  for any 
family $(a_i: i \in I)$ the infinite sum
$\su $ is defined, and  the function $\sum $ satisfies the following: 
\begin{itemize}
\item If $I= \{j,j'\}$ has two elements, then 
$\sum_{i \in I} a_i = a_j + a_{j'}$
\item If $f:I\to J$ is a bijection, and $a_i = b_{f(i)}$ for all $i\in I$, 
then $\su = \sum_{j\in J} b_j$ 
\item  Whenever $ I = \bigcup_{k\in K} J_k$
is a disjoint union of sets, and $b_k:= \sum_{j\in J_k} a_j$
for all $k\in K$, then 
$\su = \sum_{k\in K} b_k$. 
\item $x\cdot \bigl( \su\bigr) = \sum_{i\in I}(  x\cdot a_i)$, 
 $ \bigl( \su\bigr)\cdot x  = \sum_{i\in I} (a_i\cdot x)$. 
\item If all $a_i$ are $=0$, then $\su = 0$. 
\end{itemize}

\end{Definition}
\begin{Remark}
If $S$ is a complete semiring, then there are 
cardinal numbers $\lambda_1$ (the ``characteristic cardinality''
of $S$, similar to the characteristic of a field)  and 
$\lambda_S\le \max(\lambda_1, \text{cardinality of } S)$
such that
\begin{itemize}
\item
For any set $I$ there is a subset $J\subseteq I$ of cardinality $\le \lambda_1$
such that $\sum_{i\in I} 1  =\sum_{i\in J} 1 $
\item
For any set $I$ there is a subset $J\subseteq I$ of cardinality $\le \lambda_S$
such that $\su =\sum_{i\in J} a_i$
\end{itemize}
This makes it possible to view the function $\sum$ as a set rather
than a proper class. 
\end{Remark}

\begin{Definition}
A complete semiring is called d-complete 
if for all $(a_i: i=1,2,3,\ldots)$ we have: 
\begin{quote}
If $a_1 + \cdots + a_k = c$ for all $k$,
then $\sum_{i\in \{1,2,\ldots\}} a_i = c$. 
\end{quote}
Thus, in a d-complete semiring at least the ``discrete convergence'' will
be respected by countable sums. 

$\nn$ (the natural numbers together with infinity) is  a d-complete semiring. 
If $S$ is a d-complete semiring, and $\Sigma$ is any ``alphabet'', then the set
$S \langle\!\langle \Sigma^* \rangle\!\rangle$, the formal  power series 
over $\Sigma$ with coefficients in $S$, will again naturally form a 
d-complete semiring.

\end{Definition}

In \cite{Kuich:1991}, the notion of $\omega$-continuity was developed,
which is a variant of d-completeness.

\begin{Definition}

An ``ordered semiring'' is a semiring with a partial order $\le$ such that 
the operations $+$ and $\cdot$ are weakly monotone in both arguments,
and $0$ is the least element.   In particular, this implies $a\le a+x$ for all $a,x$. 

\end{Definition}

On every semiring (or even any monoid)
we can define a transitive reflexive relation  $\le$ by 
$$a \le b \ \Leftrightarrow \ \exists x (a+x=b)$$
This relation is a quasiorder, called the ``natural quasiorder'' on $(S,+)$.

\begin{Fact}
Let $(S,+,\cdot,0,1)$ be a semiring.
The following are equivalent: 
\begin{enumerate}
\item There exists an ordered semiring $T$ such that $S$ is isomorphic to 
a subsemiring of $T$. 
\item ``$S$ can be ordered'', i.e.:  there exists a partial  order on $S$ making $S$ into an ordered semiring.
\item The natural quasiorder is antisymmetric (i.e., a partial order).
\item For all $a,x,y\in S$:  ( $a+x+y=a$ implies $a+x=a$). 
\end{enumerate}
\end{Fact}

\bigskip

\begin{Definition}  A complete ordered semiring is called ``finitary'' if: 
For any $(a_i:i\in I)$, the sum $\su$ is the least  upper bound of all finite subsums: 
$$ \su = \sup \left\{ \sum_{i\in F} a_i :  F\subseteq I \mbox{ finite}\right\}$$
\end{Definition}

\begin{Fact}  Every finitary semiring is d-complete.  

The characteristic cardinality of any  finitary semiring is at most
$\aleph_0$.
\end{Fact}

\section*{Main theorem}

\begin{Theorem} \label{main}
 Let $(S,{+},{\cdot}, 0,1,{\le})$  be an ordered
semiring.  Then there is a finitary (complete ordered) semiring $\bar
S$ such that $S$ is a subsemiring of $\bar S$ (with the induced
order).    

Moreover, the construction $S \to \bar S$ is ``universal'', that is: 
For all  finitary semirings $T$ and all embeddings $f:S\to T$
(respecting the semiring operations and the order) there is a unique
embedding $\bar f: \bar S \to T$ (respecting semiring operations,
order, and therefore also $\sum$) such that $\bar f $ extends $f$. 

In other words: Let $\O$ be the category of ordered semirings (where
maps have to preserve order and the semiring operations), and let $\F$
be the category of finitary semirings, $U:\F\to \O$ the forgetful
functor. Then the map $S\mapsto \bar S$ is a left adjoint to $U$. 
\\
However, see example \ref{no}.
\end{Theorem}

\begin{proof}

Let $S^*$ be the free (multiplicative)
monoid over $S$ (consisting of all formal 
words with letters in the alphabet  $S$, including the empty word), 
and let $\nnp$ be the set of formal polynomials
over the alphabet $S$  with coefficients in $\mathbb N$
(i.e., maps from $S^*$ into $\mathbb N$ which are 
0 except on a finite set), with 
pointwise addition, and multiplication defined by the Cauchy product.

Let $\nn:= \{0,1,2,\ldots, \infty\}$, and let 
$\nnq$ be the set of formal power series over the alphabet $S$ with 
coefficients in $\nn$ (i.e.,
the set of all functions from $S^*$ to $\nn$, again with addition defined
pointwise, and the Cauchy product.

$\nnp$ is a semiring.   As as set (but not as a semiring),
$S$ is naturally embedded in $\nnp$ through a map $e:S\to \nnp$. 
There  is a natural semiring homomorphism (the ``evaluation map'') 
 $\varphi: \nnp \to S$, satisfying $\varphi(e(s))=s$  for all $s\in S$. 

$\nnq$ is a finitary complete semiring (with the natural order), and
$\nnp$ is to a subsemiring of $\nnq$. 
 We will find $\bar S$ as a homomorphic image of $\nnq$. 

$$
\begin{array}{ccc}
\nnp & \subseteq & \nnq \\
e\bigg\uparrow\ \ \varphi \bigg\downarrow && \kappa \bigg\downarrow\\
S   & \cdots\to & \bar S\\
\end{array}
$$

\newcommand{\ls}{\lesssim}
For $r,s\in \nnq$ define $ r\ls s $ iff
\begin{quote} 
For all $p\in \nnp$ with $p\le_{\nnq} r$ there exists $q\in \nnp$, $q\le_{\nnq} s$,
and $\varphi(p)\le_S \varphi(q)$.
\end{quote}

Now check that the relation $r \sim s \ :\Leftrightarrow r \ls s \ \& \ s
\ls r$ is a congruence relation on the semiring $\nnq$. 
Moreover, this congruence relation also respects infinite sums, so 
$$\bar S:= \nnq/{\sim}$$ 
is a complete semiring. 
 Clearly $S$ is embedded in $\bar S$ as an ordered semiring. 

It is easy to see that 
$\bar S$ is actually a complete finitary semiring.

\smallskip
For the ``moreover'' part:   note that a partially ordered semiring
can be made into a complete {\it finitary} semiring in at most one
way.   This fact helps to show that the construction $S\mapsto \bar S$
is universal.

\end{proof}

\begin{Corollary}\label{cor}
Let $(S,+,\cdot,0,1)$ be a semiring.  The following are equivalent:
\begin{enumerate}
\item $S$  can be ordered.
\item $S$  is a subsemiring of some d-complete semiring.
\item $S$ is a subsemiring of some complete finitary semiring. 
\end{enumerate}
\end{Corollary}
\begin{proof}

(1) $\Rightarrow$ (3) is proved in theorem \ref{main}. 

(3)  $\Rightarrow$ (2) is trivial.

(2)  $\Rightarrow$ (1):  Assume that $a+x+y=a$. Then 
$a = a+(x+y)+(x+y)+(x+y)+\cdots = (a+x) + (y+x) + (y+x)+\cdots  = a+x$. 
\end{proof}

My thesis contained also the following remark:  
{\footnotesize

\begin{Fact}
Let $(S,+,\cdot,0,1)$ be a semiring.
The following are equivalent:
\begin{enumerate}
\item $S$  is zero sum free (i.e., $x+y=0$ implies $x=y=0$.)
\item $S$  is a subsemiring of some complete semiring.
\end{enumerate}
\end{Fact}
\begin{proof}
(2) $\Rightarrow $ (1):  Similar to Corollary~\ref{cor}, (2) $\Rightarrow$ (1).

(1) $\Rightarrow $ (2):  Adjoin an element $\infty$, and declare $\su = \infty 
$ iff the set $\{i: a_i\not=0\}$  is infinite, or if for some $i$ we have $a_i=\infty$. 
\end{proof}
}

However, it is not clear that this completion will satisfy the 
infinite distributivity law.

\section*{Examples}

\begin{enumerate}
\item 
Let $S= {\mathbb N}$, then $\bar S = \nn$. 
\item 
Let $S$ be the semiring of finite subsets of some set $X$, together with
$X$ itself, with addition=union, multiplication=intersection.
  Then $\bar S$ is the full powerset of $X$. 
\item 
Let $S$ be the set of finite formal languages $L \subseteq \Sigma^*$ over 
an alphabet $\Sigma$, with  $L_1+ L_2 = L_1\cup L_2$, 
$L_1 \cdot L_2 = \{v\cdot w: v\in L_1, w\in L_2\}$.
\\
Then $\bar S$ is the set of all formal languages over $\Sigma^*$. 
\item 
Example of a complete but not d-complete semiring: 
$\{0, \text{finite}, \text{infinite}\}$, with the obvious operations.
\item
Example of a d-complete 
nonfinitary semiring:
$$\{0, \text{finite}, \text{countable}, \text{uncountable}\}$$
with the obvious operations. 
This semiring has uncountable characteristic cardinality.
\end{enumerate}

\begin{Remark}  Finitary semirings are d-complete and
have characteristic cardinality $\le \aleph_0$. 

The converse is true for finite semirings:  If $S$ is finite, d-complete
and has characteristic cardinality $\le \aleph_0$, then $S$ is finitary.

Example of a d-complete semiring with characteristic cardinality $\aleph_0$ 
which is not finitary: 
$$ \{0,1,2,\ldots, \infty-2, \infty-1, \infty \}$$
with the obvious operations, e.g.:  
$n + (\infty - k) = \infty$ if $n\ge k$, and $\infty - (k-n)$ otherwise
\end{Remark}

The following example shows that the restriction to {\em finitary} semirings
in theorem \ref{main} is reasonable. 
\begin{Example}\label{no}
There is no embedding $f$ from $\mathbb N$ into a
  complete semiring  $\mathfrak N$ such that
\begin{quote}
For all complete semirings $C$ and all semiring embeddings $g:{\mathbb N}\to 
C$ there is a complete semiring embedding $h:{\mathfrak N}\to C$ with
$h\circ f=g$. 
\end{quote}
\end{Example}
\begin{proof} Use $C$ with a large characteristic cardinality.
\end{proof}

\bibliographystyle{plain}
\bibliography{other}

\begin{thebibliography}{1}

\bibitem{Glazek:2002}
Kazimierz G{\l}azek.
\newblock {\em Semirings and Their Applications in Mathematics and Information
  Sciences}.
\newblock Kluwer, 2002.

\bibitem{Golan:1992}
Jonathan~S. Golan.
\newblock {\em The theory of semirings with applications in mathematics and
  theoretical computer science}.
\newblock Longman Scientific \& Technical, Harlow, 1992.

\bibitem{Golan:1999}
Jonathan~S. Golan.
\newblock {\em Semirings and their applications}.
\newblock Kluwer Academic Publishers, Dordrecht, 1999.
\newblock Updated and expanded version of {\it The theory of semirings, with
  applications to mathematics and theoretical computer science} 
  

\bibitem{Golan+Wang:1996}
Jonathan~S. Golan and Huaxiong Wang.
\newblock On embedding in complete semirings.
\newblock {\em Comm. Algebra}, 24(9):2945--2962, 1996.

\bibitem{Karner:1993}
Georg Karner.
\newblock On limits in complete semirings.
\newblock {\em Semigroup Forum}, 45(2):148--165, 1992.

\bibitem{Kuich:1991}
Werner Kuich.
\newblock Automata and languages generalized to $\omega$-continuous semirings.
\newblock {\em Theoret. Comput. Sci.}, 79(1, (Part A)):137--150, 1991.
\newblock Algebraic and computing treatment of noncommutative power series
  (Lille, 1988).

\bibitem{Kuich:1995}
Werner Kuich.
\newblock Representations and complete semiring morphisms.
\newblock {\em Inform. Process. Lett.}, 56(6):293--298, 1995.

\bibitem{Kuich:1997}
Werner Kuich.
\newblock Semirings and formal power series: their relevance to formal
  languages and automata.
\newblock In {\em Handbook of formal languages, Vol.\ 1}, pages 609--677.
  Springer, Berlin, 1997.

\end{thebibliography}

\bigskip
\hrule width 5cm
\bigskip
{
Martin Goldstern\\
Algebra, TU Wien\\
Wiedner Hauptstr 8--10 / 118\\
A-1040 Wien\\
\"Osterreich / Austria (Europe)\\
\
\tt goldstern@tuwien.ac.at\\
http://info.tuwien.ac.at/goldstern/}

\end{document}